\documentstyle[11pt]{article}
\begin{document}
 \baselineskip 0.22in

\title{\textbf{Convexity Conditions of Kantorovich Function and Related Semi-infinite Linear Matrix Inequalities
 }}
\author{YUN-BIN ZHAO \thanks{School of Mathematics, University of Birmingham,
Edgbaston, Birmingham B15 2TT,  United Kingdom (\tt
zhaoyy@maths.bham.ac.uk). }}

\date{}

\maketitle

 \textbf{Abstract.}
   The Kantorovich function $(x^TAx)( x^T A^{-1} x)$,
where $A$ is a positive definite matrix,   is not convex in general.
From  matrix/convex analysis point of view, it is
interesting  to address the question: When is this function convex?
  In this paper, we
investigate the convexity of this function by the condition number
of its matrix. In 2-dimensional space, we prove that the
Kantorovich function is convex if and only if the condition number
of its matrix is bounded above by $3+2\sqrt{2}, $ and thus the
convexity of the function with  two variables can be completely characterized  by the
condition number. The upper bound `$3+2\sqrt{2} $' is turned out to be a
necessary condition for the convexity of
Kantorovich functions in any finite-dimensional spaces. We also point out that when the  condition
number of the matrix (which can be any dimensional) is less than or
equal to $\sqrt{5+2\sqrt{6}}, $ the Kantorovich function is convex.
Furthermore, we prove that this general sufficient convexity condition can be remarkably improved in 3-dimensional
space. Our analysis shows that the convexity of the function is
closely related to some  modern optimization topics such as the
semi-infinite linear matrix inequality or `robust positive
semi-definiteness' of symmetric matrices.  In fact, our main result
 for 3-dimensional cases has been proved   by finding an
explicit solution range to some semi-infinite linear matrix
inequalities.\\


\newpage

\section{Introduction}
Denote by $$K(x)= (x^T A x)( x^T A^{-1} x) $$ where $x\in R^n$ and $A$
is a given $n\times n$ symmetric, positive definite real matrix with
eigenvalues $ 0<
 \lambda_1\leq \lambda_2 \leq  \cdots\leq  \lambda_n. $   Then, we have the following Kantorovich inequality (see
e.g. \cite{GR59,HJ86,  H65, K48,S59}):
$$\|x\|_2^4 \geq \frac{4\lambda_1\lambda_n}{\lambda_1^2+\lambda_n^2}
(x^T A x x^T A^{-1} x) ~\textrm{ for any }  x\in R^n. $$ This
inequality and its variants have many applications in matrix
analysis, statistics, numerical algebra, and optimization (see e.g.
\cite{BW75, D08, G88,G93, H65, HZ05, KR81, K75, S95, W87, WAS97,
Y06}). In this paper,  $K(x)$ is referred to as the `Kantorovich
function'. While $K(x)$ has been widely studied and generalized to
different forms in the literature, from matrix/convex  analysis
point of view some fundamental questions associated with this
function remain open or are not fully addressed. For instance, when
is this function convex? Is it possible to characterize its
convexity by the condition number of its matrix?

Let us first take a look at a simple example: Let $A$ be a $2\times
2$ diagonal matrix with diagonal entries $1 $ and $6.$ Then $A^{-1}$
is  diagonal with diagonal entries $1$ and $1/6$, and it is easy to
verify that in this case the function  $ K(x) $ is not convex (since
its Hessian matrix is not positive semi-definite at (1,1)). Thus,
the Kantorovich function is not convex in general. The aim of this paper is
trying to address the above-mentioned questions, and to develop some
sufficient and/or necessary  convexity conditions for this function.

We will prove that the Kantorovich function in 2-dimensional space
is convex if and only if the condition number of its matrix is less
than or equal to $3+2\sqrt{2}. $ Therefore, the convexity of this
function can be characterized completely by the condition number of
its matrix, and thus the aforementioned questions are affirmatively
and fully answered in  2-dimensional space.  For higher
dimensional cases, we prove that the upper bound `$3+2\sqrt{2}$' of
the condition number is a necessary condition for the convexity of
$K(x).$ In another word, if $K(x) $ is convex, the condition number of $A$
must not exceed this constant. On the other hand, we show that if
the condition number of the matrix is less than or equal to
$\sqrt{5+2\sqrt{6}}, $ then $K(x)$ must be convex. An immediate question is  how tight this bound of the condition number is.  Can such a general sufficient condition be improved?     A remarkable progress in this direction can be
achieved at least for  Kantorovich functions in 3-dimensional space, for which
we prove that the bound $\sqrt{5+2\sqrt{6}}$ can be improved to
$2+\sqrt{3}. $  The proof of such a result is far more than
straightforward.  It is worth mentioning that we do
not know at present whether the result in 2-dimensional space
remains valid in higher dimensional spaces. That is, it is not clear
whether or not there is a constant $\gamma$ such that the following
result holds: $K(x)$ is convex if and only if the condition number
is less than or equal to $\gamma.$ By our sufficient, and necessary
conditions, we may conclude that if such a constant $\gamma$ exists,
then $\sqrt{5+2\sqrt{6}}\leq \gamma\leq 3+2\sqrt{2}.$

The investigation of  this paper not only yields some interesting
results and new understanding for the function $K(x),$  but raises
some new challenging questions and links to certain   topics
of modern optimization and matrix analysis as well. First, the
analysis of this paper indicates  that the convexity issue of $K(x)$
is directly related to the so-called (semi-infinite) linear matrix
inequality problem which is one of the central topics of modern
convex optimization and has found broad applications in control
theory, continuous and discrete optimization, geometric distance
problems, and so on (see e.g. \cite{BEB94, D09, D08, EN00, HK93,
NN94}). In fact,  the  convexity condition of $K(x)$ can be
formulated as a  semi-infinite linear matrix inequality problem. In
 3-dimensional space, we will show that how such a semi-infinite
matrix inequality is explicitly solved in order to develop a
convexity condition for $K(x).$

The convexity issue of $K(x)$ can be also viewed as the
so-called `\emph{robust positive semi-definiteness}' of certain
symmetric matrices, arising naturally from the analysis to the Hessian
matrix of $K(x).$  The typical robust problem on the positive
semi-definiteness can be stated as follows: Let the entries of a
matrix be multi-variable functions, and some of the variables can
take any values in some intervals. The question is what range of
values the other variables should take such that the matrix is
positive semi-definite. Clearly, this question is also referred  to
as a robust optimization problem or robust feasibility/stability
problem \cite{AZ08, BN98, BN09, BS04, BS06, E97, E98, MMBBLW07}.
`Robust positive semi-definiteness' of a matrix may be stated in
different versions. For instance, suppose that the entries of the
matrix are uncertain, or cannot be given precisely, but the ranges
(e.g.  intervals) of the possible values of entries are known. Does
the matrix remain positive semi-definite when its entries vary in
these ranges?  Our analysis shows that the study of the convexity of
$K(x)$ is closely related to these topics.

Finally,  the convexity issue of $K(x)$  may stimulate the study of
more general functions than $K(x).$ Denote by $q_A(x) =\frac{1}{2}
x^T A x . $ Notice that $ q^*_A(x) =\frac{1}{2} x^T A^{-1} x =
q_{A^{-1}}(x) $ is the Legendre-Fenchel transform of $ q_A(x) $ (see
e.g. \cite{A84, BV04,   HUL01, R70}). The Kantorovich function can
be rewritten as
$$ K(x)=  4 q_A(x) q_A^* (x) = 4 q_A(x) q_{A^{-1}} (x). $$
Thus, $K(x)$ can be viewed as the product of the quadratic form
$q_A(x)$ and its Legendre-Fenchel transform $q^*_{A}(x), $ and  can be
viewed also as a special case of the product of quadratic forms. Thus,
one of the generalization of $K(x)$ is the product
$\widetilde{K}(x)= h(x) h^*(x), $ where $h$ is a convex function and
$h^*(x)$ is the Legendre-Fenchel transform of $h.$ The product of
convex functions has been exploited in the field of global
optimization under the name of multiplicative programming problems.
However, to our knowledge, the function like $ \widetilde{K}(x)$ has
not been discussed in the literature.   It is worth mentioning that the recent study for  the product of univariate convex functions and the product of quadratic forms  can be found in \cite{HU06, Z09a}, respectively.

This paper is organized as follows. In Section 2, we establish some
general sufficient, necessary conditions for the convexity of
$K(x),$ and point out  that the convexity issue of $K(x)$ can be
reformulated as a semi-infinite linear matrix inequality or robust
positive semi-definiteness of  matrices. In Section 3, we prove that
the convexity of  $K(x)$ in 2-dimensional space can be
completely characterized by the condition number of its matrix.  In
the Section 4, we prove an improved sufficient convexity condition
for the function $K(x) $ in 3-dimensional space by finding an
explicit solution range to  a class of semi-infinite linear matrix
inequalities. Conclusions are given in the last section.

\emph{Notation:} Throughout this paper, we use $A \succ 0 ~(\succeq
0) $ to denote the positive definite
  (positive semi-definite) matrix. $\kappa(A)$ denotes the
  condition number of $A$, i.e., the ratio of its largest and
  smallest eigenvalues:  $ \kappa(A) = \lambda_{\max}(A) / \lambda_{\min}
  (A).$ $q_A(x) $ denotes the quadratic form $(1/2)x^T A x.$

\section{Sufficient, necessary  conditions for the convexity of $K(x)$}  First of all, we note that a sufficient convexity condition for $K(x)$ can be obtained
by Theorem 3.2 in \cite{Z09a} which claims that the following result holds for the product of any two positive
definite quadratic forms: \emph{Let $A, B  $ be two $n\times n$
matrices and $A, B\succ 0. $ If $\kappa(B^{-1/2}AB^{-1/2}) \leq
5+2\sqrt{6}$, then the function $\left(\frac{1}{2} x^T A x\right)
\left( \frac{1}{2} x^TBx\right) $ is convex.}  By setting $B=A^{-1}$
and noting that in this case
$\kappa(B^{-1/2}AB^{-1/2})=\kappa(A^2),$ we have  the following
result.\\

\textbf{Theorem 2.1.} \emph{ Let $A$ be any $n\times n$ matrix and
$A \succ 0$. If $ \kappa(A) \leq \sqrt{5+2\sqrt{6}} ,$ then the
Kantorovich function $K(x)$  is convex.}\\

Thus, an immediate question arises: What is a necessary
condition for the convexity of $K(x)$?  This question is answered by
the next result of this section.

Let $A$ be an $n\times n$
positive definite matrix.  Denote by
\begin{equation}\label{fff} f(x): =q_A(x) q_{A^{-1}}(x) = \left( \frac{1}{2} x^T A x \right) \left(
\frac{1}{2} x^T A^{-1} x \right) =\frac{1}{4} K(x).
\end{equation}
  Clearly, the convexity of $K(x)$ is exactly
the same as that of $f(x).$ Since $f$ is twice continuously
differentiable in $R^n,$  the function $f$ is convex if and only if its Hessian
matrix $\nabla^2 f $  is positive semi-definite at any point in
$R^n.$
 It is easy to verify that the Hessian matrix of $f$
 is given by
\begin{equation}\label{DDDfff} \nabla^2 f(x)= q_A(x) A^{-1}+ q_{A^{-1}} (x) A + A x x^T
A^{-1}+ A^{-1} x x^T A.
\end{equation}
 Since $A$ is positive
definite, there exists an orthogonal matrix $U $ (i.e., $U^T U =I)$
such that \begin{equation} \label{AAAAAA}  A= U^T \Lambda U, ~~~
A^{-1} =  U^T\Lambda^{-1} U,
\end{equation}
where $\Lambda $ is a diagonal matrix whose diagonal entries
are eigenvalues of $A$ and arranged in non-decreasing order, i.e.,
$$ \Lambda=  \textrm{diag }(\lambda_1  ,  \lambda_2, ...,    \lambda_n ),  ~~~ 0<\lambda_1 \leq \lambda_2 \leq \cdots
\leq \lambda_n. $$
 By
setting $y=U x,$ we have
$$q_A(x) = \frac{1}{2} x^T(U^T \Lambda U) x
=q_\Lambda (y),  ~~~ q_{A^{-1}} (x)=\frac{1}{2} x^T(U^T \Lambda^{-1}
U) x
 = q_{\Lambda^{-1}}  (y). $$
 Notice that
$$q_\Lambda(y)= \frac{1}{2}y^T \Lambda y= \frac{1}{2} \sum_{i=1}^n \lambda_i y_i^2, ~~ q_{\Lambda^{-1}}(y)= \frac{1}{2}y^T \Lambda^{-1} y= \frac{1}{2}
\sum_{i=1}^n \frac{1}{\lambda_i} y_i^2.$$ Thus, we have
{\small
\begin{eqnarray}\label{111}  & & q_\Lambda(y) \Lambda^{-1} +
q_{\Lambda^{-1}}(y) \Lambda
 \nonumber \\
& & =   \left(
                                 \begin{array}{cccc}
                                   \frac{1}{2} \sum_{i=1}^n  \left(\frac{\lambda_i}{\lambda_1}+ \frac{\lambda_1}{\lambda_i} \right) y_i^2   &  &  &  \\
                                    & \frac{1}{2} \sum_{i=1}^n  \left(\frac{\lambda_i}{\lambda_2}+\frac{\lambda_2}{\lambda_i} \right)   y_i^2    &  &  \\
                                    &  &  \ddots  &   \\
                                     &   &   &  \frac{1}{2} \sum_{i=1}^n  \left(\frac{\lambda_i}{\lambda_n}+  \frac{\lambda_n}{\lambda_i}\right)  y_i^2   \\
                                 \end{array}
                               \right) \nonumber \\
 & & =  \left(
                                 \begin{array}{cccc}
                                  y_1^2 + \frac{1}{2} \sum_{i=2}^n  \Delta_{1i} y_i^2   &  &  &  \\
                                    & y_2^2+ \frac{1}{2} \sum_{i=1, i\not=2}^n  \Delta_{2i}  y_i^2    &  &  \\
                                    &  &  \ddots  &   \\
                                     &   &   & y_n^2+  \frac{1}{2} \sum_{i=1, i\not=n}^n  \Delta_{ni}  y_i^2   \\
                                 \end{array}
                               \right) \nonumber\\
\end{eqnarray} }
where
\begin{equation}\label{Delta-1}  \Delta_{ij}
= \frac{\lambda_j}{\lambda_i}+ \frac{\lambda_i}{\lambda_j},~~ i,j=1, ..., n, ~
i\not=j. \end{equation}
 It is evident that
\begin{equation}\label{Delta-2} \Delta_{ij}=\Delta_{ji} \textrm{ and }\Delta_{ij}\geq 2  \textrm{ for any } i
\not=j.
\end{equation}  Thus, throughout the
paper we do not distinguish  $\Delta_{ij} $ and $\Delta_{ji},$ and
 we denote by $\Delta = (\Delta_{ij})$, the
$\left(\frac{1}{2}(n-1)n\right)$-dimensional vector whose components
are $\Delta_{ij}$'s where $1\leq i< j \leq n .$  On the other hand,
it is easy to verify that
\begin{equation} \label{222}  \Lambda^{-1}
y y^T \Lambda  + \Lambda y y^T \Lambda^{-1} = \left(
                            \begin{array}{cccc}
                                    2y_1^2  &  \Delta_{12}  y_1y_2  &   \cdots &   \Delta_{1n} y_1 y_n  \\
                                     \Delta_{21} y_2y_1   &   2 y_2^2    &  \cdots &  \Delta_{2n} y_2 y_n \\
                                    \vdots & \vdots &  \ddots  &  \vdots \\
                                    \Delta_{n1} y_ny_1  &  \Delta_{n2} y_ny_2 &  \cdots  &   2 y_n^2   \\
                                 \end{array}
                               \right)
 \end{equation}
 where $\Delta_{ij} $ is given by (\ref{Delta-1}).  Therefore,  by (\ref{DDDfff}), (\ref{AAAAAA}),  (\ref{111}) and
 (\ref{222}), we have
\begin{eqnarray} \label{fffHHH} \nabla^2 f (x)  & = &  U^T  \left(q_\Lambda (y) \Lambda^{-1}+
q_{\Lambda^{-1}} (y) \Lambda + \Lambda y y^T \Lambda^{-1}+
\Lambda^{-1} y y^T \Lambda\right)  U \nonumber \\
 &   =&  U^T H_n(\Delta, y) U,
\end{eqnarray}
where
 {\small  \begin{eqnarray} \label{HHH} & &  H_n(\Delta, y) \nonumber \\
   & & = \left[
  \begin{array}{cccc}
    3 y_1^2+\sum_{ i=2}^n \frac{1}{2} \Delta_{1i}y_i^2  &  \Delta_{12}y_1y_2 &  \cdots & \Delta_{1n}y_1y_n  \\
     \Delta_{21}y_2y_1  &  3 y_2^2+ \sum_{i=1, i\not=2}^n \frac{1}{2} \Delta_{2i}y_i^2    &    \cdots
      &  \Delta_{2n}y_2y_n  \\
    \vdots & \vdots &  \ddots & \vdots \\
    \Delta_{n1}y_ny_1 &  \Delta_{n2}y_ny_2 & \cdots & 3y_n^2+ \sum_{i=1, i\not=n}^n  \frac{1}{2}
    \Delta_{ni}y_i^2 \\
  \end{array}
\right]. \nonumber \\
\end{eqnarray} }
Form (\ref{fffHHH}),  we see that  $  \nabla^2 f (x) \succeq 0$   for
any $x\in R^n $ if and only  if the vector $\Delta=(\Delta_{ij})$
satisfies that
\begin{equation}\label{HHH-NMI}  H_n(\Delta, y) \succeq 0 ~~\textrm{ for any } y\in
R^n .
 \end{equation}
Hence,    $f(x)$ is convex  (i.e., $K(x)$
is convex) if and only if (\ref{HHH-NMI}) holds. The
following observation is  useful for the proof of Theorem 2.3 of
this section.\\

 \textbf{Lemma 2.2.} \emph{Let $0< \gamma < \delta  $ be two constants. For any $t_1, t_2 \in
[\gamma, \delta],$ we have
$$ 2\leq  \frac{t_1}{t_2} + \frac{t_2}{t_1} \leq
 \frac{\gamma}{\delta} +\frac{\delta}{\gamma}. $$}

\emph{Proof.}  Let $g(\nu) =\nu+\frac{1}{\nu},$ where $ \nu =
\frac{t_1}{t_2}$. Since $t_1, t_2\in [\gamma, \delta]$, it is
evident that $\nu \in [\frac{\gamma}{\delta}, \frac{\delta}{\gamma}]
,$ in which $g(\nu) $ is convex. Clearly, the minimum value of
$g(\nu)$ attains at $\nu=1, $ i.e, $g(v)\geq 2. $   The maximum
value of $g(\nu)$ attains at one of   endpoints  of the interval
$[\frac{\gamma}{\delta}, \frac{\delta}{\gamma}]$, thus $  g(\nu)
\leq \frac{\gamma}{\delta}+\frac{\delta}{\gamma}, $ as desired.  $ ~~~~ \Box$\\

Noting that $0< \lambda_1 \leq \cdots \leq \lambda_n,$ by Lemma 2.1
and the definition of $\Delta_{ij}$, we have the  following
relation: $\Delta_{ij}\leq \Delta_{k_1k_2},\textrm{ for any }
k_1\leq i,j\leq k_2 \textrm{ and }i\not= j, $ where $1\leq k_1,
k_2\leq n.$ For instance, we have that $\Delta_{ij}\leq \Delta_{1n}
$ for any $i\not= j .$ We now prove a necessary convexity condition which is tight for general $K(x).$\\

\textbf{Theorem 2.3.} \emph{Let $A$ be any $n\times n$ matrix and $A
\succ 0$. If $K(x) $ is convex, then $\kappa (A) \leq 3+2\sqrt{2}.
$}
\\

\emph{Proof.}  Assume that $K(x) $ is convex and thus $f(x),$  given
as (\ref{fff}), is convex. It follows from (\ref{fffHHH})  and
(\ref{HHH}) that $H_n(\Delta, y)\succeq 0$ for any $y\in R^n $.
Particularly, $H_n(\Delta, \widetilde{y})\succeq 0 $ for
$\widetilde{y} = e_i+e_j , $ where $i\not=j $ and $e_i, e_j$ are the
$i$th and $j$th columns of the $n\times n$ identity matrix,
respectively. Notice that
  $$ H_n(\Delta, \widetilde{y})=  \left(
     \begin{array}{ccccccccc}
       \textbf{*} &         &                           &    &        &       &                          &      & \\
         & \ddots  &                           &    &        &       &                          &       & \\
         &         & 3+\frac{1}{2}\Delta_{ij}  &    &        &       & \Delta_{ij}              &       &  \\
         &         &                           & \textbf{*}  &        &       &                          &       & \\
         &         &                           &    & \ddots &       &                          &       &  \\
         &         &                           &    &        & \textbf{* }    &                          &       & \\
         &         &      \Delta_{ji}          &    &        &       & 3+\frac{1}{2}\Delta_{ji} &       & \\
         &         &                           &    &        &       &                          &\ddots &   \\
         &         &                           &    &        &       &                          &       & \textbf{*}  \\
     \end{array}
   \right)_{n\times n}
   $$
where all stared entries are positive numbers.
 Thus,
$H_n(\Delta, \widetilde{y}) \succeq 0$ implies that its principle
submatrix
$$  \left(
\begin{array}{cc}
3 +\frac{1}{2} \Delta_{ij}  &  \Delta_{ij} \\
 \Delta_{ji} &  3+ \frac{1}{2} \Delta_{ji} \\
 \end{array}
  \right)  $$ is positive semi-definite. Since the diagonal entries  of this  submatrix are positive,
  the submatrix is positive semi-definite if and only if
its determinant is nonnegative, i.e.,
$$ \left(3 +\frac{1}{2} \Delta_{ij} \right)\left(3+ \frac{1}{2} \Delta_{ji}\right) -
\Delta_{ij}\Delta_{ji} \geq 0,$$ which, by (\ref{Delta-2}), can be
written as $ 9 + 3\Delta_{ij} - \frac{3}{4} \Delta_{ij}^2 \geq 0. $
Thus, $ \Delta_{ij} \leq 6. $ Since $e_i$ and $e_j$ can be any
columns of the identity matrix, the inequality $\Delta_{ij}\leq 6 $
holds for any $i,j=1,..., n$ and $i\not=j.$  Lemma 2.1 implies that
$\Delta_{1n} = \max_{1\leq i,j\leq n, i\not=j} \Delta_{ij},$ and
hence
$$\kappa(A)+\frac{1}{\kappa(A)} = \frac{\lambda_n}{\lambda_1}+\frac{\lambda_1}{\lambda_n} = \Delta_{1n} = \max_{1\leq i,j\leq n, i\not=k}
\Delta_{ij} \leq 6, $$ which is equivalent to $\kappa (A)\leq
3+2\sqrt{2}. $   $~~~ \Box$

Notice  that for each fixed $y$, $H_n (\Delta,y)
\succeq 0 $ is a linear matrix inequality (LMI) in $\Delta.$   Since $y$ is any
vector in $R^n,$ the system (\ref{HHH-NMI}) is actually a
semi-infinite system of LMIs. Recall that  for a given function
$g(u,v): R^{p}\times R^q \to R$, the semi-infinite inequality (in
$u$)   is defined as $ g(u,v)\geq 0 $ for any $ v\in S, $
where $S \subseteq R^q $ is a set containing infinite many points.
Such inequalities have been widely used in robust control \cite{BEB94, EN00} and
so-called semi-infinite programming problems \cite{HK93}, and they can
be also interpreted as   robust feasibility or stability problems
when the value of $v$ is uncertain or cannot be given precisely, in
which case $S$ means all the possible values of $v.$   Recently, robust problems have wide applications in such areas as
 mathematical programming, structure
design, dynamic system, and financial optimization
\cite{AZ08,BN98,BN09,BS06, E97,E98, MMBBLW07}.
 Since   $ \Delta $ should be in
certain range such that $H_n(\Delta, y) \succeq 0 $ for any $y\in
R^n, $    (\ref{HHH-NMI}) can be also called `robust positive
semi-definiteness' of the matrix $H_n.$

Any sufficient convexity condition of $K(x)$ can provide some
explicit solution range to the system (\ref{HHH-NMI}).  In fact, by
Theorem 2.2, we immediately have the next result.\\

\textbf{Corollary  2.4.} \emph{Consider the semi-infinite linear matrix
inequality (in $ \Delta$): {\small $$\left[
  \begin{array}{cccc}
    3 y_1^2+\sum_{ i=2}^n \frac{1}{2} \Delta_{1i}y_i^2  &  \Delta_{12}y_1y_2 &  \cdots & \Delta_{1n}y_1y_n  \\
     \Delta_{21}y_2y_1  &  3 y_2^2+ \sum_{i=1, i\not=2}^n \frac{1}{2} \Delta_{2i}y_i^2    &    \cdots
      &  \Delta_{2n}y_2y_n  \\
    \vdots & \vdots &  \ddots & \vdots \\
    \Delta_{n1}y_ny_1 &  \Delta_{n2}y_ny_2 & \cdots & 3y_n^2+ \sum_{i=1, i\not=n}^n  \frac{1}{2}
    \Delta_{ni}y_i^2 \\
  \end{array}
\right] \succeq 0 $$ }
  for any   $y\in R^n, $ where
$\Delta_{ij}$'s satisfy (\ref{Delta-2}).  Then the solution set
$\Delta^*$ of the above system is nonempty, and any vector $\Delta=
(\Delta_{ij})$ where $ \Delta_{ij} \in \left[2, \sqrt{5+2\sqrt{6}}
\right]$ is
  in   $\Delta^*.$} \\

Conversely, any range of the feasible solution $\Delta $ to the
semi-infinite linear matrix inequality (\ref{HHH-NMI}) can provide a
sufficient condition for the convexity of $K(x)$. This idea is used
to prove an improved sufficient convexity condition for
 $K(x)$ in 3-dimensional space (see Section 4 for details).\\

Combining Theorems 2.2 and 2.3 leads to the following corollary.\\

\textbf{Corollary 2.5.} \emph{Assume that there exists a constant,
denoted by $\gamma^*$, such that the following statement is true:
 $K(x)$ is
convex if and only if $\kappa(A)\leq \gamma^*. $  Then such a
constant must satisfy that
$ \sqrt{5+2\sqrt{6}} \leq \gamma^* \leq 3+2\sqrt{2}.$}\\

However, the question is: Does such a constant exist? If the answer
is `yes', we  obtain a complete characterization of the convexity of
$K(x)$ by merely the condition number. In the next section, we prove
that this question can be fully addressed for 2-dimensional
Kantorovich functions, to which the constant is given by $ \gamma^*=
3+2\sqrt{2} $ (as a result,   the bound given in Theorem 2.3 is
tight). For higher dimensional cases, the answer to this question is
not clear at present. However, for 3-dimensional cases we can prove
 that the lower bound $\sqrt{5+2\sqrt{6}} $ can be significantly improved (see Section 4 for details).

\section{Convexity characterization  for  $K(x)$ in 2-dimensional space}

We now consider the  Kantorovich function with two variables
and prove that the necessary condition in Theorem 2.3 is also
sufficient for this case.\\

\textbf{Theorem 3.1.}  \emph{Let $A $ be any $2\times 2 $ positive
definite matrix. Then $ K(x)   $ is convex if and only if $\kappa
(A) \leq 3+2\sqrt{2}.$}
\\

 \emph{Proof.}  Notice
that there exists an orthogonal matrix  $U$ (i.e., $U^T U=I$) such
that
$$ A=  U^T\left[\begin{array}{cc}
\beta_1 & 0 \\
0 & \beta_2
\end{array} \right] U,~~  A^{-1}= U^T\left[\begin{array}{cc}
\frac{1}{\beta_1} & 0 \\
0 &\frac{1}{ \beta_2}
\end{array}\right] U , $$
where
 $\beta_1, \beta_2$ are  eigenvalues
 of $A$. By using the same notation of Section 2, and setting $y=
 Ux
 $  and  $n=2$ in (\ref{HHH}), we have
\begin{equation}\label{HHH222} H_2(\Delta, y)   = \left[\begin{array}{cc} 3 y_1^2+   \frac{1}{2} \Delta_{12}y_2^2   &  \Delta_{12}y_1y_2   \\
     \Delta_{21}y_2y_1  &  3 y_2^2+   \frac{1}{2} \Delta_{21}y_1^2
  \end{array}
\right].
\end{equation}
For 2-dimensional cases, the vector $\Delta $ is reduced to the
scalar $\Delta = \Delta_{12} . $  If $ f(x) $ given by (\ref{fff})
is convex in $R^n, $   it follows from (\ref{fffHHH}) that
$H_2(\Delta, y)\succeq 0 $ for any   $y\in R^2.$ In particular, it
must be positive semi-definite at $y= e =(1, 1)$, thus
$$ H _2(\Delta, e) =    \left[\begin{array}{cc}
3+  \frac{1}{2} \Delta_{12}  &
  \Delta_{12}   \\
  \Delta_{21}
 & 3+  \frac{1}{2} \Delta_{21}
\end{array}\right] \succeq 0.
$$
Since $\Delta_{12}=\Delta_{21},$ it follows that $ \det  H_2(\Delta,
e) =   \left(3+ \frac{1}{2}  \Delta_{12}\right)^2 - \Delta_{12} ^2
\geq 0 , $ i.e.
\begin{equation}\label{betabeta} 9+3
\Delta_{12} -\frac{3}{4} \Delta_{12} ^2 \geq 0.
\end{equation}
Conversely, if (\ref{betabeta}) holds, we can prove that $f$ is
convex.  Indeed,  we see that  the diagonal entries of $H_2(\Delta,
y)$ are positive,  and that
\begin{eqnarray*} \det H_2(\Delta, y) & = & \left(3y_1^2+  \frac{1}{2}
\Delta_{12} y_2^2\right) \left( 3 y_2^2+  \frac{1}{2} \Delta_{12}
y_1^2\right) - \Delta_{12}^2
y_1^2y_2^2  \\
& = & 9 y_1^2 y_2^2+ \frac{3}{2} \Delta_{12} (y_1^4+y_2^4)
-\frac{3}{4} \Delta_{12}^2
y_1^2y_2^2  \\
& \geq  &  9 y_1^2 y_2^2+ 3\Delta_{12} y_1^2y_2^2  -\frac{3}{4}
\Delta_{12}^2
y_1^2y_2^2  \\
& = &   \left( 9+3 \Delta_{12} -\frac{3}{4}
\Delta_{12} ^2 \right) y_1^2y_2^2
  \geq       0.
\end{eqnarray*}
The first inequality above follows from the fact $y_1^4+y_2^4\geq 2
y_1^2y_2^2 $ and the second inequality follows from
(\ref{betabeta}). Thus, $ H_2(\Delta, y) \succeq 0 $ for any $y\in
R^2, $ which  implies that $f$ is convex.

Therefore,   $f(x)$ is convex  (i.e., $K(x)$ is convex) if and only
if (\ref{betabeta}) holds. Notice that the roots of the quadratic
function $9+3t- \frac{3}{4} t^2 =0 $ are  $ t^*_1= -2 $ and $
t^*_2=6 . $  Notice that $\Delta_{12} \geq 2 $ (see (\ref{Delta-2})).  We conclude that
(\ref{betabeta}) holds if and only if $ \Delta_{12} \leq 6. $ By the
definition of $\Delta_{12},$ we have
$$ \Delta_{12} = \frac{\beta_1}{\beta_2}+\frac{\beta_2}{\beta_1}
=\kappa (A)+\frac{1}{\kappa(A)},
$$  Thus the inequality $\Delta_{12}\leq 6 $  can be written as  $\kappa (A)^2-6\kappa(A)
+1\leq 0,$ which is equivalent to $  \kappa (A)     \leq 3+2\sqrt{2}.
$  $ ~~~~ \Box $\\

It is worth stressing that the above result can be also obtained by
solving the semi-infinite linear matrix inequality (\ref{HHH-NMI}).
In fact, by Theorem 2.3, it suffices to prove that $\kappa(A) \leq
3+2\sqrt{2}$ is  sufficient for the convexity of   $K(x)  $ in 2-dimensional space.   Suppose that  $\kappa(A)\leq 3+2\sqrt{2},$ which is equivalent
to $\Delta_{12}=\frac{\beta_1}{\beta_2}+\frac{\beta_2}{\beta_1} \leq
6. $  Thus, $\Delta=\Delta_{12} \in  [2,6] . $ We now prove that $K
(x) $ is convex.  Define
$$ \upsilon (\Delta_{12}, y) := \det H_2(\Delta, y). $$
   Differentiating the function
with respect to $\Delta_{12},$ we have
$$ \frac{\partial \upsilon (\Delta_{12}, y) } {\partial \Delta_{12}}
= \det \left[\begin{array}{cc}   \frac{1}{2}  y_2^2   &   y_1y_2   \\
     \Delta_{21}y_2y_1  &  3 y_2^2+   \frac{1}{2} \Delta_{21}y_1^2
  \end{array}
\right] + \det \left[\begin{array}{cc} 3 y_1^2+   \frac{1}{2} \Delta_{12}y_2^2   &  \Delta_{12}y_1y_2   \\
     y_2y_1  &     \frac{1}{2}  y_1^2
  \end{array}
\right] .$$ Differentiating it again, we have
\begin{eqnarray*}
\frac{\partial^2 \upsilon (\Delta_{12}, y) } {\partial^2
\Delta_{12}}
= \det \left[\begin{array}{cc}   \frac{1}{2}  y_2^2   &   y_1y_2   \\
    y_2y_1  &      \frac{1}{2} y_1^2
  \end{array}
\right] + \det \left[\begin{array}{cc}    \frac{1}{2}  y_2^2   &  y_1y_2   \\
     y_2y_1  &     \frac{1}{2}  y_1^2
  \end{array}
\right] = -\frac{3}{2} y_1^2y_2^2 \leq 0.
\end{eqnarray*}
Therefore, for any given $y\in R^2$, the function $\upsilon
(\Delta_{12}, y)$ is concave with respect to $\Delta_{12},$ and
hence the minimum value of the function attains at one of  the
endpoints of the interval $ [2,6],$ i.e., \begin{equation}
\label{vvv-ineq} \upsilon (\Delta_{12}, y) \geq \min\{\upsilon (2,
y), \upsilon (6, y)\} . \end{equation}  Notice that
$$ \upsilon (2, y) =\det  \left[\begin{array}{cc} 3 y_1^2+    y_2^2   &  2y_1y_2   \\
    2y_2y_1  &  3 y_2^2+   y_1^2
  \end{array}
\right] $$ where the matrix is diagonally dominant and its diagonal
entries are nonnegative for any given $y\in R^2,$  and thus it is
positive semi-definite. This implies that $\upsilon (2, y) \geq 0 $
for any $y \in R^2. $  Similarly,   $$ \upsilon (6,
y) = \det \left[\begin{array}{cc} 3 y_1^2+    3y_2^2   &  6y_1y_2   \\
    6y_2y_1  &  3 y_2^2+   3y_1^2
  \end{array}
\right] \geq 0,  $$ since the matrix is diagonally dominant the
diagonal entries are nonnegative for any $y\in R^2. $   Therefore,
it follows from (\ref{vvv-ineq}) that $ \upsilon (\Delta_{12}, y)
\geq 0 $ for any $ \Delta_{12}\in [2,6] \textrm{ and } y\in R^2. $
This implies that when $\Delta_{12}\in [2,6] $ the matrix
$H_2(\Delta, y) \succeq 0 $ for any $y\in R^2, $  and hence $K(x) $
is convex.

\section{An improved convexity condition for  $K(x)$ in 3-dimensional space}
In this section, we prove that for  $K(x)$ in 3-dimensional space, the upper
bound of the  condition number given in Theorem 2.2 can be improved
to $  \leq 2+\sqrt{3}. $ To prove this result, we try to
find an explicit solution range  to a class of semi-infinite linear
matrix inequalities, which  will immediately yield an improved
sufficient convexity condition  for  $K(x) $ in 3-dimensional space. First, we
give some useful inequalities.\\

\textbf{Lemma 4.1.} \emph{If $ 2 \leq   \delta_{1}, \delta_2,
\delta_3 \leq 4,  $ then all the functions below are nonnegative:}
\begin{eqnarray} \chi_1(\delta_1, \delta_2, \delta_3)  & : = & ~ 6\delta_3+  \delta_1\delta_2-\frac{1}{2}\delta_3\delta^2_{2}
\geq 0, \label{inequality-01}\\
   \chi_2 (\delta_1, \delta_2, \delta_3) & : = &  ~6\delta_1 + \delta_3 \delta_2-\frac{1}{2}
\delta_1 \delta^2_{2}  \geq 0,\label{inequality-02} \\
 \chi_3 (\delta_1, \delta_2, \delta_3) & : =  & ~ 6 \delta_{2} +  \delta_1  \delta_3 - \frac{1}{2}
\delta_2 \delta_1^2   \geq 0, \label{inequality-03}\\
  \chi_4 (\delta_1, \delta_2, \delta_3) & : = &  ~6\delta_3+
\delta_1\delta_2-\frac{1}{2} \delta_3\delta_1^2 \geq 0, \label{inequality-04} \\
   \psi( \delta_1, \delta_2, \delta_3) & : =& 12 + \delta_1\delta_2 \delta_3
-\delta_1 ^2-\delta_2^2-\delta_3^2 \geq 0.\label{inequality-05}
\end{eqnarray}

\emph{Proof.} For any given $\delta_1, \delta_3\in [2,4],$ we
consider the quadratic function (in $t$): $  \rho(t) = 6 \delta_3
+ \delta_1 t -\frac{1}{2} \delta_3 t^2  $   which is concave in $t.$
Let $t_2^*$ be the largest root of $\rho(t) =0.$ Then
$$t_2^*= \frac{  \delta_1 + \sqrt{ \delta_1^2+12
\delta_3^2} } {\delta_3} =  \frac{\delta_1}{\delta_3} + \sqrt{
\left(\frac{\delta_1}{\delta_3} \right)^2 + 12 }  \geq \frac{2}{4} +
\sqrt{ \left(\frac{2}{4}\right)^2+ 12} = 4, $$  where the inequality
follows from the fact that $\frac{\delta_1}{\delta_3} \geq
\frac{2}{4}.$  It is easy to check that the least root of $\rho
(t) =0 $ is non-positive, i.e., $t^*_1 \leq 0.$ Thus, the interval
$[2,4]\subset [t_1^*, t_2^*] $ in which the quadratic function
$\rho (t) \geq 0.$  Since $2\leq \delta_2\leq 4, $ we conclude
that
$$ \chi_1(\delta_1, \delta_2, \delta_3) =  \rho(\delta_2) =  6\delta_3 + \delta_1
\delta_2-\frac{1}{2} \delta_3 \delta^2_2  \geq 0. $$ Thus
(\ref{inequality-01}) holds.   All inequalities
(\ref{inequality-02})-(\ref{inequality-04}) can be proved by the
same way, or simply by
  exchanging
the role of $\delta_1, \delta_2 $  and  $\delta_3$ in
(\ref{inequality-01}).

 We now prove (\ref{inequality-05}).
It is easy to see that $\psi$ is concave with respect to its every
variable. The minimum value of   $\psi$ with respect to $\delta_1$
attains at the boundary of the interval $[2,4].$ Thus, we have
\begin{equation}\label{DDDD-01} \psi (\delta_1, \delta_2,
\delta_3)   \geq    \min\{\psi ( 2, \delta_2,\delta_3), \psi (4,
\delta_2, \delta_3 )\}  ~ \textrm{ for any }\delta_1, \delta_2,
\delta_3\in [2,4].
\end{equation}
Notice that
$ \psi ( 2, \delta_2,\delta_3) $ and $ \psi (4, \delta_2, \delta_3 ) $   are
concave with respect to $\delta_2\in[2,4].$ Thus we have that
\begin{eqnarray}  \psi ( 2, \delta_2,\delta_3)  & \geq &  \min\{\psi ( 2,
2,\delta_3), \psi ( 2, 4,\delta_3)\}~~\textrm{  for any }
\delta_2,\delta_3\in [2,4],  \label{DDDD-02} \\
 \psi (4, \delta_2,
\delta_3 ) & \geq  & \min\{\psi ( 4, 2,\delta_3), \psi ( 4,
4,\delta_3)\}~~\textrm{  for any } \delta_2,\delta_3\in [2,4].
\label{DDDD-03}
\end{eqnarray}
Similarly,   $ \psi ( 2, 2,\delta_3), \psi ( 2, 4,\delta_3), \psi (
4, 2,\delta_3) $ and $\psi ( 4, 4,\delta_3) $ are concave in
$\delta_3. $ Thus,
\begin{eqnarray*} \psi ( 2, 2,\delta_3) \geq
\min\{\psi ( 2, 2,2), \psi ( 2, 2,4 \} &  = &  \min\{8,  4)\} >0 ~
\textrm{ for any }  \delta_3 \in [2,4], \\
  \psi ( 2, 4,\delta_3) \geq
\min\{\psi ( 2, 4,2), \psi ( 2, 4,4 \} & = &  \min\{4,  8)\} >0 ~
\textrm{ for any }  \delta_3 \in [2,4],  \\
  \psi ( 4, 2,\delta_3) \geq
\min\{\psi ( 4, 2,2), \psi ( 4, 2,4 \} & =&  \min\{4,  8)\} >0  ~
\textrm{ for any }  \delta_3 \in [2,4],\\
  \psi ( 4, 4,\delta_3)   \geq
\min\{\psi ( 4, 4, 2), \psi ( 4, 4,4 \} & = &  \min\{8,  8)\} >0  ~
\textrm{ for any }  \delta_3 \in [2,4]. \end{eqnarray*}
 Thus, combining
(\ref{DDDD-01})-(\ref{DDDD-03}) and the last four
inequalities above yields (\ref{inequality-05}).   $ ~~~~ \Box$\\

We now focus on developing explicit solution range  to certain
semi-infinite linear matrix inequalities which will be used later to
establish an improved sufficient convexity condition for
 $K(x) $ in 3-dimensional space.

\subsection{The solution range to semi-infinite linear matrix
inequalities}

Consider the following $3\times 3$ matrix whose entries are the
functions  in $ (\omega,\alpha, \beta): $
  \begin{equation}\label{MMMM}
M(\omega, \alpha, \beta)=\left[
  \begin{array}{ccc}
    3 + \frac{1}{2} \omega_1\alpha^2 + \frac{1}{2} \omega_2\beta^2 &  \omega_1\alpha &  \omega_2\beta  \\
     \omega_1\alpha  &   \frac{1}{2} \omega_1 + 3 \alpha^2 + \frac{1}{2}
     \omega_3\beta^2
      &  \omega_3\alpha \beta  \\
    \omega_2\beta &  \omega_3\alpha\beta &  \frac{1}{2}
    \omega_2 +\frac{1}{2} \omega_3\alpha^2 +3 \beta^2\\
  \end{array}
\right] \\ \end{equation}    where  $\omega
=(\omega_1, \omega_2, \omega_3) \in R^3, $
$\omega_{j}  \in [2,6] $ for $j=1,2,3,$ and  $\alpha, \beta\in [-1, 1]. $   We are interested in
finding the range of $ \omega_1, \omega_2, $ and $\omega_3 $ in
the interval  $[2,6] $  such that
$$M(\omega, \alpha, \beta)
\succeq 0 \textrm{ for any } \alpha, \beta \in [-1, 1], $$ which is
a semi-infinite linear matrix inequality. To this end, we seek the
condition for $\omega, $ under which all the principle minors of $M$
are nonnegative for any $\alpha, \beta \in [-1,1]. $

 First of all, we see from (\ref{MMMM}) that all diagonal entries (which are the first order principle minors) of $M$  are
positive in the intervals considered. Secondly, since $ \omega_j \in
[2,6] (j=1,2,3),$ and $\alpha \in [-1,1],$ it is easy to see that
 \begin{eqnarray*} 3 + \frac{1}{2} \omega_1\alpha^2  & = & \frac{1}{2} (6+\omega_1 \alpha^2)   \geq   \frac{1}{2} (\omega_1+ \omega_1\alpha^2) \geq
\omega_1 |\alpha |, \\
 \frac{1}{2} \omega_1 + 3 \alpha^2 & =&
\frac{1}{2}(\omega_1+6 \alpha^2) \geq
\frac{1}{2}(\omega_1+\omega_1 \alpha^2) \geq \omega_1
|\alpha |.
\end{eqnarray*}
Therefore,  the second order principle submatrix  of (\ref{MMMM})
$$ \left[
  \begin{array}{cc}
    3 + \frac{1}{2} \omega_1\alpha^2 + \frac{1}{2} \omega_2\beta^2 &  \omega_1\alpha    \\
     \omega_1\alpha  &   \frac{1}{2} \omega_1 + 3 \alpha^2 + \frac{1}{2}
     \omega_3\beta^2
  \end{array}
\right]
$$  is diagonally dominant for any  $\omega_{j}\in [2,6] (j=1,2,3)$ and $\alpha, \beta\in
[-1,1].$ Similarly, the second order principle submatrices
$$ \left[
  \begin{array}{cc}
    3 + \frac{1}{2} \omega_1\alpha^2 + \frac{1}{2} \omega_2\beta^2 &  \omega_2\beta    \\
     \omega_2\beta  &   \frac{1}{2} \omega_2  + \frac{1}{2}
     \omega_3\alpha^2+ 3 \beta^2
  \end{array}
\right],
$$ $$ \left[
  \begin{array}{cc}
    \frac{1}{2} \omega_1+3\alpha^2 + \frac{1}{2} \omega_3\beta^2 &  \omega_3\alpha \beta   \\
     \omega_3\alpha\beta  &   \frac{1}{2} \omega_2   + \frac{1}{2}
     \omega_3\alpha^2+ 3 \beta^2
  \end{array}
\right]
$$
are also diagonally dominant for any $\omega_{j}\in [2,6] (j=1,2,3) $ and
$\alpha, \beta \in [-1,1]. $ Thus, all second order principle minors
of $M$ are nonnegative in the intervals considered. It is sufficient
to find the range of $\omega=(\omega_1, \omega_2, \omega_3)
$ such  that the third order principle minor is nonnegative, i.e.,
$$\det M(\omega, \alpha, \beta) \geq 0 ~\textrm{  for any } \alpha, \beta  \in [-1,
1]. $$ However, this is not straightforward.  In order to find such
a range, we calculate up to the sixth order partial derivative of
$\det M(\omega, \alpha, \beta)$ with respect to $\alpha.$ The first
order partial derivative is given as
\begin{equation}\label{DDD111} ~~~~~ \frac{\partial \det
 M(\omega, \alpha, \beta)}{\partial\alpha}= \det M_1(\omega,\alpha,
\beta)+\det M_2( \omega, \alpha, \beta) +\det M_3( \omega, \alpha,
\beta),
\end{equation}
 where $M_i(\omega, \alpha, \beta)$ is the matrix that
coincides with the matrix $M(\omega, \alpha, \beta)$ except that
every entry in the $i$th row is differentiated with respect to
$\alpha, $ i.e.,
{\small
\begin{equation}\label{M1} ~~~~~M_1( \omega, \alpha, \beta) =
\left[\begin{array}{ccc}
                           \omega_1 \alpha & \omega_1  &   0  \\
                            \omega_1 \alpha  &  \frac{1}{2}\omega_1+3\alpha^2+\frac{1}{2} \omega_3 \beta^2 &
                            \omega_3\alpha\beta \\
                           \omega_2\beta & \omega_3 \alpha \beta &  \frac{1}{2}\omega_2+\frac{1}{2} \omega_3
                           \alpha^2 +3\beta^2\\
                         \end{array}
                       \right],
\end{equation} }
{\small
\begin{equation} \label{M2} ~~~~M_2(\omega,\alpha, \beta) = \left[
                         \begin{array}{ccc}
                           3+\frac{1}{2} \omega_1 \alpha^2+\frac{1}{2} \omega_2 \beta^2  & \omega_1 \alpha  &   \omega_2\beta  \\
                            \omega_1  &  6 \alpha  &
                            \omega_3\beta \\
                           \omega_2\beta & \omega_3 \alpha \beta &  \frac{1}{2}\omega_2+\frac{1}{2} \omega_3
                           \alpha^2 +3\beta^2\\
                         \end{array}
                       \right],
\end{equation} }
{\small
\begin{equation} \label{M3}  ~~~~ M_3(\omega,\alpha, \beta) = \left[
                         \begin{array}{ccc}
                           3+\frac{1}{2} \omega_1 \alpha^2+\frac{1}{2} \omega_2 \beta^2  & \omega_1 \alpha  &   \omega_2\beta   \\
                             \omega_1 \alpha  &  \frac{1}{2}\omega_1+3\alpha^2+\frac{1}{2} \omega_3 \beta^2 &
                            \omega_3\alpha\beta  \\
                           0 &  \omega_3  \beta  &   \omega_3
                           \alpha  \\
                         \end{array}
                       \right].
\end{equation} }
 Similarly,  the notation $M_{ij}(\omega, \alpha, \beta) $ means the resulting matrix by differentiating, with respect to $\alpha,$ every entry
 of the $i$th and $j$th rows of  $M(\omega, \alpha, \beta), $  respectively;  $M_{ijk} (\omega, \alpha, \beta) $ means the
matrix obtained by differentiating every entry of the $i$th, $j$th
and $k$th rows of $M(\omega, \alpha, \beta)$, respectively. All
other matrices $M_{ijkl...} (\omega, \alpha, \beta)$ are understood
this way.
 In particular, the matrix such as $M_{iij} (\omega, \alpha, \beta)  $ with some identical indices means the  matrix obtained by
 differentiating every entry of  the $i$th row of $M (\omega, \alpha, \beta)$ twice and differentiating every entry of its $j$th
 row once.  Notice that every entry of  $M(\omega,
\alpha, \beta)$ is at most quadratic in $\alpha. $ Thus if we
differentiate a row  three times or more, the resulting matrix contains
a row with all entries zero, and hence its determinant is equal to
zero. For example, $\det M_{1113}(\omega, \alpha, \beta) =0. $

  Clearly, the second order partial derivative is given as
follows:
\begin{eqnarray}\label{DDD222} & &  \frac{\partial^2 \det M(\omega, \alpha, \beta)}{\partial^2 \alpha} \nonumber  \\&
& =   \frac{\partial\det M_1(\omega, \alpha, \beta)}{\partial
\alpha} + \frac{\partial \det M_2(\omega, \alpha, \beta)}{\partial
\alpha} +
\frac{\partial \det M_3(\omega, \alpha, \beta) }{\partial \alpha}, \nonumber\\
& & = \sum_{j=1}^3 \det M_{1j} ( \omega, \alpha, \beta) +
\sum_{j=1}^3 \det
M_{2j} (\omega, \alpha, \beta) + \sum_{j=1}^3 \det M_{3j} (\omega, \alpha, \beta) \nonumber\\
& &=  \sum_{i=1}^3\sum_{j=1}^3 \det M_{ij} (\omega, \alpha, \beta)  \nonumber\\
& & = \det M_{11}(\omega, \alpha, \beta)  + \det M_{22} (\omega,
\alpha, \beta) + \det M_{33} (\omega, \alpha, \beta)  \nonumber \\
& & ~+2\det M_{12}
(\omega, \alpha, \beta) +2\det M_{13} (\omega, \alpha, \beta) +2\det M_{23} (\omega,
\alpha, \beta),
\end{eqnarray}
where the last equality follows from the fact $M_{ij} (\omega,
\alpha, \beta) =M_{ji} (\omega, \alpha, \beta) .$ By differentiating
(\ref{DDD222})   and noting that $\det M_{111}=\det M_{222}=\det
M_{333}=0,$ we have
\begin{eqnarray}\label{DDD333}
   & &  \frac{\partial^3 \det M(\omega, \alpha, \beta)}{\partial^3
\alpha} =    \sum_{i=1}^3\sum_{j=1}^3 \sum_{k=1}^3\det M_{ijk}
(\omega, \alpha, \beta) \nonumber \\
&   & =    3\det M_{112} (\omega, \alpha, \beta) +3\det M_{113}
(\omega, \alpha, \beta) +3\det M_{122} (\omega, \alpha, \beta)
    \nonumber \\& & ~~ +6\det M_{123} (\omega, \alpha, \beta) +3\det
M_{133}(\omega, \alpha, \beta) +3\det M_{223}(\omega, \alpha, \beta) \nonumber \\
& & ~~ +3\det M_{233}(\omega, \alpha, \beta).
\end{eqnarray}
By differentiating it again and noting that $\det M_{iiij} =0, $ we
have
\begin{eqnarray} \label{DDD444}
  & &  \frac{\partial^4 \det M(\omega,\alpha,
\beta)}{\partial^4 \alpha}    = \sum_{i=1}^3\sum_{j=1}^3
\sum_{k=1}^3 \sum_{l=1}^3 \det M_{ijkl} (\omega, \alpha, \beta), \nonumber \\
&   &  =  6\det M_{1122}(\omega, \alpha, \beta) +  6\det
M_{1133}(\omega,
\alpha, \beta) + 6\det M_{2233} (\omega, \alpha, \beta) \nonumber \\
 & & ~~~  + 12
\det M_{1123}(\omega, \alpha, \beta)    +  12 \det M_{1223} (\omega,
\alpha, \beta) + 12 \det M_{1233} (\omega, \alpha, \beta).
\end{eqnarray}
Finally,
\begin{eqnarray}
  \frac{\partial^5 \det M(\omega,\alpha,
\beta)}{\partial^5 \alpha}   & =  &   \sum_{i=1}^3\sum_{j=1}^3
\sum_{k=1}^3 \sum_{l=1}^3 \sum_{p=1}^3 \det M_{ijklp}
(\omega,\alpha, \beta), \label{DDD555}\\
  \frac{\partial^6 \det M(\omega,\alpha,
\beta)}{\partial^6 \alpha}   & = &   \sum_{i=1}^3\sum_{j=1}^3
\sum_{k=1}^3 \sum_{l=1}^3 \sum_{p=1}^3\sum_{q=1}^3 \det M_{ijklpq}
(\omega, \alpha, \beta). \label{DDD666}
\end{eqnarray}
Our first technical result is given as follows.\\

 \textbf{Lemma 4.2.} \emph{Let $ \omega_1, \omega_2,\omega_3 \in [2, 4] $ and $\alpha, \beta \in [-1, 1].$ Then the  function $ \frac{\partial^4
\det M(\omega, \alpha, \beta)}{\partial^4 \alpha} $ is convex with
respect to $\alpha, $ and $ \frac{\partial^4 \det M(\omega, \alpha,
\beta)}{\partial^4 \alpha} \geq 0 . $}\\

 \emph{Proof.} Let $\beta\in [-1,1] $ and $\omega_1, \omega_3, \omega_2 \in [2,4] $ be arbitrarily given. Define $g(\alpha) = \frac{\partial^4
\det M(\omega, \alpha, \beta)}{\partial^4 \alpha}, $ which is a
function in $ \alpha.$  As we mentioned earlier, if we differentiate
a row of $M (\omega, \alpha, \beta) $ three times or more, the
determinant of the resulting matrix  is equal to zero.
 Thus, the nonzero terms on the right-hand side of (\ref{DDD666}) are only those
matrices obtained by differentiating every row of $M (\omega,
\alpha, \beta)$  exactly twice, i.e, the terms
$$\det M_{ijklpq} (\omega, \alpha, \beta) =\det M_{112233}(\omega, \alpha, \beta)
 = \det \left[
 \begin{array}{ccc}
 \omega_1 & 0 & 0 \\
 0 & 6  & 0 \\
 0 & 0 & \omega_3 \\
 \end{array}
  \right]=6\omega_1\omega_3
> 0, $$ and hence  $ g''(\alpha) = \frac{\partial^6 \det
M(\omega, \alpha, \beta)}{\partial^6 \alpha} > 0. $ This implies
that  $g(\alpha) =  \frac{\partial^4 \det M(\omega, \alpha,
\beta)}{\partial^4 \alpha} $ is convex with respect to $\alpha \in
[-1,1].$ Notice that $g'(\alpha) = \frac{\partial^5 \det M(\omega,
\alpha, \beta)}{\partial^5 \alpha}. $ We now prove that $$ g'(0) =
\frac{\partial^5 \det M(\omega, \alpha, \beta)}{\partial^5
\alpha}\Big{|}_{\alpha=0} =0.$$
  When the matrix $M $ is differentiated 5 times with respect to $\alpha,$ there are only three possible cases
  in which we   have  nonzero
determinants for the resulting matrices.

Case 1:  Rows 1, 2 are differentiated twice, and row 3 once. In this case we have
$$ \det M_{11223} (\omega, \alpha, \beta) = \det \left[
  \begin{array}{ccc}
     \omega_1  &   0 &  0   \\
     0   &  6
      &  0   \\
    0  &  \omega_3 \beta &    \omega_3\alpha \\
  \end{array}
\right] = 6 \omega_1\omega_3\alpha. $$

Case 2:  Rows 1, 3 are differentiated twice, and row 2 once. We have
$$ \det M_{11233} (\omega, \alpha, \beta) =  \det \left[
  \begin{array}{ccc}
     \omega_1  &   0 &  0   \\
     \omega_1   &  6\alpha
      &  \omega_3\beta  \\
    0  &   0 &    \omega_3\\
  \end{array}
\right] = 6 \omega_1\omega_3\alpha.  $$

Case 3: Rows 2, 3 are done twice, and row 1 once. Then
$$\det  M_{12233} (\omega, \alpha, \beta) =  \det  \left[
  \begin{array}{ccc}
     \omega_1\alpha  &  \omega_1  0 &  0   \\
     0  &  6 &  0   \\
    0  &   0 &    \omega_3\\
  \end{array}
\right] = 6 \omega_1\omega_3\alpha . $$  Clearly,
\begin{eqnarray*} g'(\alpha)  & = &  m_1 \det M_{11223}(\omega, \alpha, \beta) + m_2 \det M_{11233}(\omega, \alpha, \beta) +
m_3 \det M_{12233} (\omega, \alpha, \beta) \\
& = &  6(m_1+m_2+m_3)\omega_1\omega_3 \alpha, \end{eqnarray*}
where $m_1, m_2, m_3 $ are  positive integers due to the duplication
of the terms in (\ref{DDD555}), such as $M_{11223} (\omega, \alpha,
\beta) =M_{12123} (\omega, \alpha, \beta) =M_{22131}(\omega, \alpha,
\beta) .$ Therefore, $g'(0)=0. $ By the convexity of $g(\alpha),$
the  minimum value of $g$ attains  at $\alpha=0.$ We now prove that this
minimum value is nonnegative, and hence $g(\alpha)\geq 0.$ Indeed,
it is easy to see that
 $$  M_{1223} (\omega, \alpha, \beta) = \left[
                 \begin{array}{ccc}
                   \omega_1\alpha &  \omega_1 &  0 \\
                  0  & 6  & 0  \\
                   0  & \omega_3\beta & \omega_3\alpha \\
                 \end{array}
               \right], ~M_{1123}(\omega, \alpha, \beta) = \left[
                 \begin{array}{ccc}
                   \omega_1 &  0 &  0 \\
                   \omega_1 & 6\alpha & \omega_3\beta \\
                   0  & \omega_3\beta & \omega_3\alpha \\
                 \end{array}
               \right]
$$
$$  M_{1233}(\omega, \alpha, \beta)  = \left[
                 \begin{array}{ccc}
                   \omega_1\alpha &  \omega_1 &  0  \\
                   \omega_1 & 6\alpha & \omega_3\beta \\
                   0  & 0 & \omega_3  \\
                 \end{array}
               \right] , $$
               $$  M_{1133} (\omega, \alpha, \beta)  = \left[
                 \begin{array}{ccc}
                   \omega_1  &   0  &   0  \\
                   \omega_1\alpha  &   \frac{1}{2} \omega_1+ 3\alpha^2 + \frac{1}{2} \omega_3 \beta^2 & \omega_3\alpha\beta  \\
                   0  & 0 & \omega_3  \\
                 \end{array}
               \right] ,$$
$$  M_{1122}  (\omega, \alpha, \beta) = \left[
                 \begin{array}{ccc}
                   \omega_1  & 0  &  0 \\
                  0  & 6  & 0  \\
                   \omega_2\beta  & \omega_3\alpha \beta & \frac{1}{2} \omega_2+\frac{1}{2} \omega_3\alpha^2+3\beta^2 \\
                 \end{array}
               \right],
$$

$$  M_{2233} (\omega, \alpha, \beta) = \left[
                 \begin{array}{ccc}
                   3+\frac{1}{2} \omega_1\alpha^2 + \frac{1}{2} \omega_2\beta^2 & \omega_1\alpha & \omega_2\beta \\
                  0  & 6  & 0  \\
                   0   & 0  &   \omega_3  \\
                 \end{array}
               \right].
$$
%Notice that
%  $$ \det M_{1223}
%(\omega, \alpha, \beta) \Big{|}_{\alpha=0} = 0, ~ \det M_{1133} (\omega, \alpha, \beta) \Big{|}_{\alpha=0} =
% \frac{1}{2}(\omega_1^2
% \omega_3+ \omega_1 \omega^2_{3}\beta^2),  $$
%$$\det M_{1123} (\omega,\alpha, \beta) \Big{|}_{\alpha=0} = - \omega_1 \omega_3^2  \beta^2,
% ~ \det M_{1122} (\alpha, \beta) \Big{|}_{\alpha=0} =
% 3\omega_1\omega_2 +18 \omega_1\beta^2, $$
%$$  \det M_{1233} (\alpha, \beta) \Big{|}_{\alpha=0} =  -\omega_1^2
% \omega_3,
%  ~~\det M_{2233} (\omega, \alpha, \beta) \Big{|}_{\alpha=0} =
% 18 \omega_3  +3 \omega_2 \omega_3 \beta^2. $$
 Therefore, by (\ref{DDD444}) we have
  \begin{eqnarray*}  & &   g(0)
   =   \frac{\partial^4 \det M(\omega, \alpha, \beta)}{\partial^4
\alpha}\Big{|}_{\alpha=0} \\
 &   & =    12\left[\det M_{1123} (\omega, \alpha, \beta)|_{\alpha=0} + \det
M_{1223}(\omega,
\alpha, \beta)|_{\alpha=0}  +  \det M_{1233} (\omega, \alpha, \beta)|_{\alpha=0}\right] \\
 & & ~~ +6\left[\det
M_{1122} (\omega, \alpha, \beta) |_{\alpha=0} +  \det
M_{1133}(\omega, \alpha,
\beta)|_{\alpha=0}   + \det M_{2233} (\omega, \alpha, \beta)|_{\alpha=0}\right] \\
&  & =  18 \left(6\omega_1 +  \omega_2
\omega_3-\frac{1}{2} \omega_1 \omega^2_{3} \right)\beta^2 +18
 \left(6\omega_3 +
 \omega_1 \omega_2-\frac{1}{2} \omega_1^2 \omega_3
\right).
\end{eqnarray*}
By Lemma 4.1, when $\omega_1, \omega_2,\omega_3\in [2,4],$
we have \begin{eqnarray*} 6\omega_1 + \omega_2 \omega_3-\frac{1}{2}
\omega_1 \omega^2_{3} & =&  \chi_2(\omega_1,
\omega_3,\omega_2 )  \geq 0, \\
 6\omega_3 + \omega_1
\omega_2-\frac{1}{2} \omega_1^2 \omega_3 & = &
\chi_3(\omega_1, \omega_3,\omega_2 ) ( = \chi_4(\omega_1, \omega_2,\omega_3 )) \geq 0 . \end{eqnarray*}
 Therefore,
$g(\alpha)\geq g(0) \geq 0 $ for any $\alpha \in [-1,1].$ Since
$\beta \in [-1,1],$ and $\omega_1, \omega_3, \omega_2 \in
[2,4]$ are arbitrarily given points, the desired results follows.   $ ~~~~ \Box$\\

 \textbf{Lemma 4.3.} \emph{Let $\omega_1,\omega_2,\omega_3 \in [2,4]$  and $\alpha, \beta \in [-1,1]. $ Then
$ \frac{\partial^2 \det M(\omega, \alpha, \beta)}{\partial^2 \alpha}
$ is convex with respect to $\alpha, $  and $\alpha=0 $ is a
minimizer of it, i.e. \begin{equation}\label{IIIEEE}
\frac{\partial^2 \det M(\omega, \alpha, \beta)}{\partial^2 \alpha}
\geq \frac{\partial^2 \det M(\omega, \alpha, \beta)}{\partial^2
\alpha}\Big{|}_{\alpha=0} \end{equation} for any $\alpha, \beta,
\omega $ in the above-mentioned intervals.}\\

\emph{ Proof.} The convexity of  $\frac{\partial^2 \det M(\omega,
\alpha, \beta)}{\partial^2 \alpha} $ with respect to $\alpha$ is an
immediate consequence of Lemma 4.2 since  $\frac{\partial^4 \det
M(\omega, \alpha, \beta)}{\partial^4 \alpha}$, the second order
derivative of $\frac{\partial^2 \det M(\omega, \alpha,
\beta)}{\partial^2 \alpha}$,  is nonnegative.  We now prove
(\ref{IIIEEE}).
  To this end, we first show that
\begin{equation}\label{D333} \frac{\partial^3 \det M(\omega, \alpha,
\beta)}{\partial^3 \alpha}\Big{|}_{\alpha=0} =0,
\end{equation}
which implies that $\alpha= 0$ is a minimizer of the second order
partial derivative. In fact, by (\ref{DDD333}), to calculate the third order
partial derivative with respect to $\alpha$  we
need to calculate the determinant of $M_{ijk}(\omega, \alpha,
\beta).$ Clearly,
 \begin{eqnarray*}   M_{123}(\omega, \alpha, \beta)    =     \left[
  \begin{array}{ccc}
     \omega_1\alpha  &  \omega_1   &  0   \\
     \omega_1  &  6\alpha &   \omega_3\beta   \\
    0  &  \omega_3\beta &    \omega_3 \alpha \\
  \end{array}
\right].
\end{eqnarray*}
 and thus $
\det M_{123}(\omega, \alpha, \beta) |_{\alpha=0} = 0.  $ Similarly,
we have
 \begin{eqnarray*}  M_{112}(\omega,\alpha,  \beta)  & = & \left[
  \begin{array}{ccc}
      \omega_1  &  0  &  0   \\
     \omega_1  &    6 \alpha
      &  \omega_3 \beta  \\
  \omega_2\beta &  \omega_3\alpha\beta &  \frac{1}{2}
    \omega_2 +\frac{1}{2} \omega_3\alpha^2 +3 \beta^2   \\
  \end{array}
\right],\\
   M_{113}(\omega,\alpha,  \beta) &= & \left[
  \begin{array}{ccc}
     \omega_1   &  0  &  0   \\
     \omega_1\alpha  &   \frac{1}{2} \omega_1 + 3 \alpha^2 + \frac{1}{2}
     \omega_3\beta^2
      &  \omega_3\alpha \beta  \\
   0  &  \omega_3 \beta &  \omega_3\alpha \\
  \end{array}
\right], \\
  M_{221}(\omega, \alpha,  \beta)  & = &  \left[
  \begin{array}{ccc}
     \omega_1 \alpha  &  \omega_1  &  0   \\
     0  &   6   &  0   \\
    \omega_2\beta &  \omega_3\alpha\beta &  \frac{1}{2}
    \omega_2 +\frac{1}{2} \omega_3\alpha^2 +3 \beta^2\\
  \end{array}
\right],
 \\
  M_{223}(\omega, \alpha,  \beta) & = &  \left[
  \begin{array}{ccc}
    3 + \frac{1}{2} \omega_1\alpha^2 + \frac{1}{2} \omega_2\beta^2 &  \omega_1\alpha &  \omega_2\beta  \\
      0  &   6   &  0  \\
    0  &  \omega_3 \beta &   \omega_3\alpha \\
  \end{array}
\right],  \end{eqnarray*}  which imply that  $ \det M_{112}(\omega, \alpha,
\beta) |_{\alpha=0}=  \det M_{113}(\omega, \alpha, \beta)
|_{\alpha=0}= \det M_{221}(\omega,
\alpha, \beta) |_{\alpha=0}= \det M_{223}(\omega, \alpha,  \beta)
|_{\alpha=0} =0. $  Finally, we have
 $$  M_{331}(\omega,\alpha,  \beta) =  \left[
  \begin{array}{ccc}
     \omega_1\alpha  &  \omega_1  &  0  \\
     \omega_1\alpha  &   \frac{1}{2} \omega_1 + 3 \alpha^2 + \frac{1}{2}
     \omega_3\beta^2
      &  \omega_3\alpha \beta  \\
    0  &  0  &    \omega_3 \\
  \end{array}
\right],$$
$$  M_{332}(\omega, \alpha,  \beta) =   \left[
  \begin{array}{ccc}
    3 + \frac{1}{2} \omega_1\alpha^2 + \frac{1}{2} \omega_2\beta^2 &  \omega_1\alpha &  \omega_2\beta  \\
     \omega_1   &  6 \alpha
      &  \omega_3 \beta  \\
     0  &  0  &    \omega_3
  \end{array}
\right].
$$
Clearly, we also have that  $ \det M_{331}(\omega, \alpha,  \beta) |_{\alpha=0} = \det
M_{332}(\omega,  \alpha,  \beta)|_{\alpha=0}  =0. $  From the above
calculation, by (\ref{DDD333}) we see that (\ref{D333}) holds.
 This means that $\alpha=0$ is the minimizer of $\frac{\partial^2
 M(\omega, \alpha, \beta)}{\partial^2 \alpha} $ for any given
 $\omega_1, \omega_2, \omega_3 \in[2,4]$ and $\beta \in [-1,1],$  and thus (\ref{IIIEEE}) holds.  $ ~~~~ \Box$\\

Based on the above fact, we may further show that the second order
partial derivative is positive in the underlying intervals.\\

 \textbf{Lemma 4.4.} \emph{If $\omega_1, \omega_3, \omega_2 \in [2, 4 ] $ and $\alpha, \beta\in [-1,1], $ then
$ \frac{\partial^2
 M(\omega, \alpha, \beta)}{\partial^2 \alpha} \geq 0. $}\\

\emph{Proof.}  By Lemma 4.3, $\alpha =0$ is a minimizer of $
\frac{\partial^2
 M(\omega, \alpha, \beta)}{\partial^2 \alpha}. $  Thus, it is sufficient to show that at $\alpha=0$ the second order partial derivative is nonnegative
 for any given $\beta \in [-1,1]$ and $ \omega_1,
\omega_3, \omega_2 \in [2, 4 ] . $   Indeed,
\begin{eqnarray*} & & \frac{\partial^2 \det M(\omega, \alpha, \beta)}{\partial^2 \alpha} \\
& & =  \det M_{11}(\omega, \alpha, \beta)+2[ \det
M_{12}(\omega,\alpha, \beta)  + \det M_{13}(\omega, \alpha, \beta) +
\det
M_{23}(\omega, \alpha, \beta)]\\
 & & ~~~ + \det M_{22}(\omega, \alpha,
\beta)  + \det M_{33}(\omega, \alpha, \beta)\\
&  & = \det \left[
  \begin{array}{ccc}
      \omega_1  &  0  &  0  \\
     \omega_1\alpha  &   \frac{1}{2} \omega_1 + 3 \alpha^2 + \frac{1}{2}
     \omega_3\beta^2
      &  \omega_3\alpha \beta  \\
    \omega_2\beta &  \omega_3\alpha\beta &  \frac{1}{2}
    \omega_2 +\frac{1}{2} \omega_3\alpha^2 +3 \beta^2\\
  \end{array}
\right] \\
& & ~~~ + 2 \det  \left[
  \begin{array}{ccc}
      \omega_1\alpha   &  \omega_1  &  0  \\
     \omega_1   &   6\alpha
      &   \omega_3  \beta  \\
    \omega_2\beta &  \omega_3\alpha\beta &  \frac{1}{2}
    \omega_2 +\frac{1}{2} \omega_3\alpha^2 +3 \beta^2\\
  \end{array}
\right]  \\
& & ~~~  + 2 \det  \left[
  \begin{array}{ccc}
      \omega_1\alpha  &  \omega_1  & 0   \\
     \omega_1\alpha  &   \frac{1}{2} \omega_1 + 3 \alpha^2 + \frac{1}{2}
     \omega_3\beta^2
      &  \omega_3\alpha \beta  \\
    0  &  \omega_3 \beta &    \omega_3\alpha \\
  \end{array}
\right] \\
& & ~~~  + 2\det \left[
  \begin{array}{ccc}
    3 + \frac{1}{2} \omega_1\alpha^2 + \frac{1}{2} \omega_2\beta^2 &  \omega_1\alpha &  \omega_2\beta  \\
     \omega_1   &  6  \alpha
      &  \omega_3 \beta  \\
   0  &  \omega_3 \beta &  \omega_3\alpha \\
  \end{array}
\right] \\
& & ~~~ + \det \left[
  \begin{array}{ccc}
    3 + \frac{1}{2} \omega_1\alpha^2 + \frac{1}{2} \omega_2\beta^2 &  \omega_1\alpha &  \omega_2\beta  \\
     0  &   6
      &  0  \\
    \omega_2\beta &  \omega_3\alpha\beta &  \frac{1}{2}
    \omega_2 +\frac{1}{2} \omega_3\alpha^2 +3 \beta^2\\
  \end{array}
\right] \\
&  &  ~~~ +  \det  \left[
  \begin{array}{ccc}
    3 + \frac{1}{2} \omega_1\alpha^2 + \frac{1}{2} \omega_2\beta^2 &  \omega_1\alpha &  \omega_2\beta  \\
     \omega_1\alpha  &   \frac{1}{2} \omega_1 + 3 \alpha^2 + \frac{1}{2}
     \omega_3\beta^2
      &  \omega_3\alpha \beta  \\
   0  &  0 &    \omega_3 \\
  \end{array}
\right],
\end{eqnarray*}
and hence, by setting $\alpha=0$ and by a simple calculation we have
  \begin{eqnarray*}      \frac{\partial^2 \det M(\omega, \alpha, \beta)}{\partial^2 \alpha}\Big{|}_{\alpha=0}
&  = &   9\omega_2+\frac{3}{2}
\omega_1\omega_3-\frac{3}{4}\omega_1^2 \omega_2  +
\frac{9}{2} \left(12+ \omega_1\omega_2\omega_3
-\omega_1^2-\omega_3^2-\omega_2^2 \right) \beta^2 \\
&  & ~~  +
\left(9\omega_2+\frac{3}{2}\omega_1\omega_3-\frac{3}{4}\omega_2\omega^2_{3}
\right)\beta^4.
\end{eqnarray*}
By Lemma 4.1, for any $\omega_1, \omega_2, \omega_3 \in
[2,4] ,$  we have the following inequalities:
\begin{eqnarray*} 9\omega_2+\frac{3}{2}
\omega_1\omega_3-\frac{3}{4}\omega_1^2 \omega_2  & = &
\frac{3}{2} \chi_4( \omega_1, \omega_3,\omega_2) \geq 0,\\
  9\omega_2 +
\frac{3}{2}\omega_1\omega_3-\frac{3}{4}\omega_2\omega^2_{3}
& = & \frac{3}{2}\chi_1( \omega_1, \omega_3,\omega_2) \geq 0,
 \\
    12+
\omega_1\omega_2\omega_3
-\omega_1^2-\omega_3^2-\omega_2^2 & = & \psi(\omega_1,
\omega_3,\omega_2) \geq 0 .
\end{eqnarray*}
 Thus,  the minimum value of
$\frac{\partial^2 \det M(\alpha, \beta)}{\partial^2 \alpha} $ is
nonnegative.     $ ~~~~ \Box$\\

We now prove the main result of this subsection.\\

\textbf{Theorem 4.5.} \emph{If $ \omega_1, \omega_3,\omega_2\in [2,
4]$, then    $ M(\omega, \alpha, \beta) \succeq 0  $ for any
$\alpha,\beta \in [-1,1].$}\\

\emph{Proof.} By Lemma 4.4, $ \det M(\omega, \alpha, \beta) $  is
convex with respect to $\alpha $ since its second partial derivative
is nonnegative.  We now prove the $\alpha =0$ is a minimizer of
$\det M (\omega, \alpha, \beta)$ for an arbitrarily given $\beta\in
[-1, 1] $ and $\omega_{ij}$ in $ [2,4]. $ Indeed, from (\ref{M1}),
(\ref{M2}) and (\ref{M3}), it is easy to see that
$$\det M_1(\omega, \alpha, \beta)|_{\alpha=0}  = \det M_2(\omega, \alpha, \beta)|_{\alpha=0} =
\det M_3(\omega, \alpha, \beta)|_{\alpha=0} =0,$$ and hence by
(\ref{DDD111}), we have
$$ \frac{\partial \det  M(\omega, \alpha, \beta)}{\partial \alpha} \Big{|}_{\alpha=0}
=0,$$ which together with the convexity of $\det M(\omega, \alpha,
\beta)$ implies that the minimum value of $\det M(\omega, \alpha, \beta)$ attains
 at $\alpha =0.$ Substituting $\alpha=0$ into $ \det
M(\omega, \alpha, \beta) $ we have \begin{eqnarray}\label{DET} \det
M(\omega, \alpha, \beta) \Big{|}_{\alpha =0} &  = &  \det \left[
  \begin{array}{ccc}
    3 +   \frac{1}{2} \omega_2\beta^2 &  0  &  \omega_2\beta  \\
     0   &   \frac{1}{2} \omega_1  + \frac{1}{2}
     \omega_3\beta^2
      & 0   \\
    \omega_2\beta & 0  &  \frac{1}{2}
    \omega_2   +3 \beta^2\\
  \end{array}
\right] \nonumber \\
& = & \frac{3}{2} \left(  \omega_1+ \omega_3 \beta^2\right)
\left(\frac{1}{2} \omega_2+\left(3-\frac{1}{4}
\omega_2^2\right)\beta^2 +\frac{1}{2} \omega_2\beta^4\right).
\end{eqnarray}
Notice that the quadratic function $9 -\frac{5}{2}t+ \frac{1}{16}t^2
\leq 0 $ for any $t\in  [4, 36].$ Since $\omega_2\in [2,4],$ we
have that  $\omega_2^2 \in [4, 16] \subset [4,36]$  we conclude
that
$$ \left(3-\frac{1}{4} \omega_2^2\right)^2 - 4\left(\frac{1}{2} \omega_2 \right) \left(\frac{1}{2} \omega_2\right)
= 9 -\frac{5}{2}\omega_2^2+ \frac{1}{16}\omega_2^4 \leq 0,$$
which means the determinant of the quadratic function $ \frac{1}{2}
 \omega_2+\left(3-\frac{1}{4} \omega_2^2\right) t +\frac{1}{2}
\omega_2 t^2 $ (in $t$) is non-positive, and thus $$ \frac{1}{2}
\omega_2+\left(3-\frac{1}{4} \omega_2^2\right)\beta^2
+\frac{1}{2} \omega_2\beta^4 \geq 0 ,$$ which together with
(\ref{DET}) implies  that
$$ \det M(\omega, \alpha, \beta) \geq \det M(\omega, \alpha,
\beta)\Big{|}_{\alpha=0} \geq 0. $$ Thus,  the third order principle
minor of $M (\omega, \alpha, \beta) $ is nonnegative. As we mentioned at the beginning of this subsection, under our conditions the diagonal entries of $M$  are positive, and
all  second order principle minors of $M$ are also nonnegative.
Thus,   $M(\omega, \alpha, \beta) \succeq 0 $ under our
conditions.  $ ~~~~ \Box$ \\

By symmetry, we also have the following result:\\

\textbf{Theorem 4.6.} \emph{If $\omega_1, \omega_3,\omega_2 \in
[2, 4],$ then $ P(\omega, \alpha, \beta) \succeq 0 $ and
  $  Q(\omega, \alpha, \beta)  \succeq  0 $   for any  $ \alpha,
\beta\in [-1, 1], $  where} {\small $$ P(\omega, \alpha, \beta) =
\left[
  \begin{array}{ccc}
    3\alpha^2 + \frac{1}{2} \omega_1  + \frac{1}{2} \omega_2\beta^2 &  \omega_1\alpha &  \omega_2\alpha\beta  \\
     \omega_1\alpha  &   \frac{1}{2} \omega_1\alpha^2 + 3   + \frac{1}{2}
     \omega_3\beta^2
      &  \omega_3 \beta  \\
    \omega_2\alpha \beta &  \omega_3 \beta &  \frac{1}{2}
    \omega_2\alpha^2 +\frac{1}{2} \omega_3 +3 \beta^2\\
  \end{array}
\right] ,$$ } {\small $$ Q(\omega, \alpha, \beta)=  \left[
  \begin{array}{ccc}
    3 \alpha^2+ \frac{1}{2} \omega_1\beta^2 + \frac{1}{2} \omega_2 &  \omega_1\alpha\beta &  \omega_2\alpha \\
     \omega_1\alpha \beta &   \frac{1}{2} \omega_1\alpha^2 + 3 \beta^2 + \frac{1}{2}
     \omega_3
      &  \omega_3 \beta  \\
    \omega_2\alpha &  \omega_3 \beta &  \frac{1}{2}
    \omega_2\alpha^2 +\frac{1}{2}  \omega_3 \beta^2 +3 \\
  \end{array}
\right]. $$ }\\

  \emph{Proof.}  This result can be proved by the same way   of
Theorem 4.5. However,  repeating the whole similar proof
 is too tedious. We may prove the result by symmetry.  In fact,
notice that all the analysis and results in this Subsection depend
only on the following conditions: $\omega_1, \omega_2, \omega_3\in
[2,4] $ and $\alpha, \beta\in [-1,1]. $  By permuting the rows 1 and
2, and columns 1 and 2 of $ P(\omega,\alpha, \beta),$ and by setting
the substitutions $ \omega_1= \widetilde{\omega}_{1},$ $\omega_2 =
\widetilde{\omega}_{3}$ and $  \omega_3 = \widetilde{\omega}_{2}, $
then $P(\omega,\alpha, \beta)$ is transformed to
$M(\widetilde{\omega}, \alpha, \beta) $ where $\widetilde{\omega}
=(\widetilde{\omega}_{1}, \widetilde{\omega}_{2},
\widetilde{\omega}_{3}) .$  Clearly, we have $\widetilde{\omega}_{j}
\in [2,4] ~ (j=1,2,3).$  Thus, the positive semi-definiteness of $P$
immediately follows from that of $M.$ Similarly, by swapping rows 2
and 3 and swapping columns 2 and 3 of $Q$ , and setting the
substitutions $ \omega_1 = \widetilde{\omega}_{2} ,$  $  \omega_2 =
\widetilde{\omega}_{1},$ and $ \omega_3 = \widetilde{\omega}_{3}, $
it is easy to see that $ Q(\omega, \alpha, \beta)$ can be
transformed into $P(\widetilde{\omega}, \alpha, \beta).$ Therefore,
the result follows from Theorem 4.5 immediately. $ ~~~~ \Box$

\subsection{The improved sufficient convexity condition}
For general Kantorovich functions,   Theorem 2.3 claims that when
$K(x)$ is convex, it must satisfy that $\kappa(A)\leq 3+2\sqrt{2}$
which is the necessary condition for $K(x)$ to be convex. This
necessary condition is equivalent to $2\leq \Delta_{ij} \leq 6 $ for
any $i\not=j.$  Thus any sufficient condition for the convexity of
$K(x)$ must fall into this range. Theorem 2.2 shows that if
$\kappa(A)\leq \sqrt{5+2\sqrt{6}} \approx 3.14626$ (which is
equivalent to $2 \leq \Delta_{ij} \leq
\frac{6+5\sqrt{6}}{\sqrt{5+2\sqrt{6}}} \approx 3.4641$ for any
$i\not =j $), then $K(x) $ is convex.   Based on the result of
Subsection 4.1, we  now prove that this  sufficient convexity
condition can be improved in 3-dimensional space.

 Notice that in 3-dimensional space, by (\ref{fffHHH}) and (\ref{HHH}), the positive semi-definiteness of the Hessian
matrix  of $K(x)$ is equivalent to that of
{\small    \begin{eqnarray}\label{HHH333} & &  H_3(y) \nonumber \\ & & = \left[
  \begin{array}{ccc}
    3 y_1^2+ \frac{1}{2} \Delta_{12}y_2^2 + \frac{1}{2} \Delta_{13}y_3^2 &  \Delta_{12}y_1y_2 &  \Delta_{13}y_1y_3  \\
     \Delta_{21}y_2y_1  &   \frac{1}{2} \Delta_{21}y_1^2 + 3 y_2^2 + \frac{1}{2}
     \Delta_{23}y_3^2
      &  \Delta_{23}y_2y_3  \\
    \Delta_{31}y_3y_1 &  \Delta_{32}y_3y_2 &  \frac{1}{2}
    \Delta_{31}y_1^2 +\frac{1}{2} \Delta_{32}y_2^2 +3 y_3^2\\
  \end{array}
\right]. \nonumber \\
\end{eqnarray}   }
The main result of this section is given below.\\

\textbf{Theorem 4.7.} \emph{If $\kappa(A)\leq 2+\sqrt{3}, $  then $
H_3(y) \succeq 0 $ for any $y\in R^n,$ and hence the  function
$K(x)=x^TAx x^T A^{-1} x $ is convex.}\\

\emph{ Proof.} Let $H_3(y)$ be given by (\ref{HHH333}). Clearly,
$H_3(0) \succeq 0.$ In what follows  we assume that $y=(y_1, y_2,
y_3) \not=0.$ There are three possible cases:

\emph{Case 1.}  $y_1 $ is the component with the largest absolute
value: $|y_1| \geq \max\{|y_2|, |y_3|\}. $ Denote by
$\alpha=\frac{y_2}{y_1}, ~~ \beta = \frac{y_3}{y_1}. $  Notice that
$\Delta_{ij}=\Delta_{ji}$ for any $i\not=j$ (see (\ref{Delta-2})).
By setting $\omega =\Delta=(\Delta_{12}, \Delta_{13}, \Delta_{23}),
$ it is easy to see from (\ref{HHH333}) that {\small
\begin{eqnarray*}
 H_3(y)
& = & y_1^2\left[
  \begin{array}{ccc}
    3 + \frac{1}{2} \Delta_{12}\alpha^2 + \frac{1}{2} \Delta_{13}\beta^2 &  \Delta_{12}\alpha &  \Delta_{13}\beta  \\
     \Delta_{21}\alpha  &   \frac{1}{2} \Delta_{21} + 3 \alpha^2 + \frac{1}{2}
     \Delta_{23}\beta^2
      &  \Delta_{23}\alpha \beta  \\
    \Delta_{31}\beta &  \Delta_{32}\alpha\beta &  \frac{1}{2}
    \Delta_{31} +\frac{1}{2} \Delta_{32}\alpha^2 +3 \beta^2\\
  \end{array}
\right]\\
& = & y_1^2 M (\Delta, \alpha, \beta)
\end{eqnarray*} }
where $M(\cdot, \cdot,\cdot )$ is defined by (\ref{MMMM}). Thus if
$M(\Delta, \alpha, \beta) \succeq 0 $ for any $\alpha, \beta\in
[-1,1],$ then $H_3(y) \succeq 0 $ for any $y$ such that $|y_1| \geq
\max\{|y_2|, |y_3|\}.$

\emph{Case 2.}   $y_2 $ is the component with the largest absolute
value:  $|y_2| \geq \max\{|y_1|, |y_3|.$ Denote by
$\alpha=\frac{y_1}{y_2}, ~~ \beta = \frac{y_3}{y_2}. $ Then,  we
have {\small \begin{eqnarray*}
 H_3(y)
& = & y_2^2  \left[
  \begin{array}{ccc}
    3\alpha^2 + \frac{1}{2} \Delta_{12}  + \frac{1}{2} \Delta_{13}\beta^2 &  \Delta_{12}\alpha &  \Delta_{13}\alpha\beta  \\
     \Delta_{21}\alpha  &   \frac{1}{2} \Delta_{21}\alpha^2 + 3   + \frac{1}{2}
     \Delta_{23}\beta^2
      &  \Delta_{23} \beta  \\
    \Delta_{31}\alpha \beta &  \Delta_{32} \beta &  \frac{1}{2}
    \Delta_{31}\alpha^2 +\frac{1}{2} \Delta_{32} +3 \beta^2\\
  \end{array}
\right]        \\
& = & y_2^2  P(\Delta, \alpha, \beta)
\end{eqnarray*}  } where $P(\cdot, \cdot,\cdot)$ is defined as in Theorem 4.6.
Thus  if $P(\Delta, \alpha, \beta) \succeq 0 $ for any $\alpha,
\beta\in [-1,1],$  then $H_3(y) \succeq 0 $ for any $y$ such that
$|y_2| \geq \max\{|y_1|, |y_3|\}. $

\emph{Case 3.}   $y_3 $ is the component with the largest absolute
value: $|y_3|  \geq \max\{|y_1|, |y_2|\}.$  Denote by $
\alpha=\frac{y_1}{y_3}, ~~ \beta = \frac{y_2}{y_3}. $ Then,
{\small \begin{eqnarray*}
 H_3(y)
& = & y_3^2    \left[
  \begin{array}{ccc}
    3 \alpha^2+ \frac{1}{2} \Delta_{12}\beta^2 + \frac{1}{2} \Delta_{13} &  \Delta_{12}\alpha\beta &  \Delta_{13}\alpha \\
     \Delta_{21}\alpha \beta &   \frac{1}{2} \Delta_{21}\alpha^2 + 3 \beta^2 + \frac{1}{2}
     \Delta_{23}
      &  \Delta_{23} \beta  \\
    \Delta_{31}\alpha &  \Delta_{32} \beta &  \frac{1}{2}
    \Delta_{31}\alpha^2 +\frac{1}{2}  \Delta_{32} \beta^2 +3 \\
  \end{array}
\right]    \\
& = & y_3^2 Q(\Delta, \alpha, \beta)
\end{eqnarray*}  } where $Q(\cdot, \cdot,\cdot)$ is defined as in Theorem 4.6.
Therefore, if $Q(\Delta, \alpha, \beta) \succeq 0 $ for any $\alpha,
\beta\in [-1,1],$  then $H_3(y) \succeq 0 $ for any $y$ such that
$|y_3| \geq \max\{|y_1|, |y_2|\}. $

If $\kappa(A)\leq 2+\sqrt{3} $ which is equivalent to $\Delta_{ij}\in [2,4]$ where
$i,j=1,2,3$ and $i\not=j$, by setting $\omega =\Delta,$ i.e., $(\omega_1, \omega_2,\omega_3) = (\Delta_{12}, \Delta_{13}, \Delta_{23}),$  from Theorems 4.5 and 4.6  we have that
$$ M (\Delta, \alpha, \beta) \succeq 0,  ~ P(\Delta, \alpha, \beta) \succeq 0,
 ~ Q(\Delta, \alpha, \beta) \succeq 0 \textrm{ for any }  \alpha, \beta \in[-1,1] . $$
Thus, $H_3(y) \succeq 0$ for any $y\in R^n,$ and hence $K(x)$ is
convex.  $ ~~~~ \Box$\\

The above theorem shows that the  upper bound
 of the condition number  in
Theorem 2.2 is improved to
 $ 2+\sqrt{3} ~(\approx 3.73205) $ in 3-dimensional space.  As
a  result, Corollary 2.4 can be  improved accordingly in
3-dimensional space.\\

\textbf{Corollary 4.8.} Let $n=3. $  Assume that there exists a
constant $\gamma^*$ such that the following statement is true:
$\kappa(A) \leq \gamma^*$ if and only if $K(x)$ is convex. Then the
constant $\gamma^*$  must satisfy that
$$ 2+\sqrt{3} \leq  \gamma^*  \leq  3+2\sqrt{2}. $$\\

\section{Conclusions}  The purpose of this paper is to characterize the
convexity of the Kantorovich function through only the condition
number of its matrix.  We have shown that if the Kantorovich
function is convex, the condition number of its matrix must be less
than or equal to $3+2\sqrt{2}. $  It turns out that such a necessary
condition is also sufficient for any Kantorovich functions in
2-dimensional space. For higher dimensional cases ($n\geq 3$), we
point out that a sufficient condition for $K(x)$ to be convex is
that the condition number is less than or equal to
$\sqrt{5+2\sqrt{6}}.$ Via certain semi-infinite linear matrix
inequalities, we have proved that this general sufficient convexity
condition  can be   improved to $2+\sqrt{3}$ in 3-dimensional space.
Our analysis shows that the convexity issue of the Kantorovich
function is closely related to a class of semi-infinite linear
matrix inequality problems, and it is also related to some robust
feasibility/semi-definiteness problems for certain matrices. Some
interesting and challenging questions on $K(x) $ are worthwhile for
the future work.  For instance, can we obtain a complete convexity
characterization  for $n\geq 3? $ In another word, for $ n\geq 3, $
does there exist a constant $\gamma^*$ such that $ \kappa(A) \leq
\gamma^* $ if and only if $K(x)$ is convex? If the answer to these
questions is 'yes', what would be the explicit value or formula
(which might depend on $n$) for this constant?

\end{document}